\documentclass[10.5pt,reqno]{amsart}

\usepackage{marginnote,amssymb,mathrsfs,latexsym,verbatim,pdfsync,color,mathabx,tikz}

\usepackage [utf8]{inputenc}
\usepackage[colorlinks=true, linkcolor=blue]{hyperref}

\def\bbC{{\mathbb C}}

\def\bbN{{\mathbb N}}

\def\bbR{{\mathbb R}}

\def\bbZ{{\mathbb Z}}

\def\cA{{\mathcal A}}

\def\cD{{\mathcal D}}

\def\cH{{\mathcal H}}

\def\cL{{\mathcal L}}

\def\cN{{\mathcal N}}

\def\cS{{\mathcal S}}
\def\cT{{\mathcal T}}

\def\cW{{\mathcal W}}

\def\cZ{{\mathcal Z}}

\def\Re{\operatorname{Re}}%{\, \text{\rm Re\,}}
\def\Im{\operatorname{Im}}%{\, \text{\rm Im\,}}

\def\la{\langle}
\def\ra{\rangle}

\def\eps{\varepsilon}
\def\l{\lambda} 
 
\def\vp{\varphi}
\def\ov{\overline}
\def\p{\partial}

\def\ms{\medskip}

\def\endpf{\medskip\hfill $\Box$

\smallskip %\ms

}  
\def\epf{\endpf}

\def\Hol{\operatorname{Hol}}

\def\har{\operatorname{har}}

\def\supp{\operatorname{supp}}

\def\sgn{\operatorname{sgn}}

\def\itb{\item[{\tiny $\bullet$}]}

\def\Hol{\operatorname{Hol}}

\def\CB{\color{black} }
\def\CR{\color{red} }

\newtheorem{thm}{Theorem}[section]
\newtheorem{prop}[thm]{Proposition}

\newtheorem{lem}[thm]{Lemma}
\newtheorem{defn}[thm]{Definition}
\newtheorem{remark}[thm]{Remark}

\usepackage{graphicx}

\usepackage{amsmath}

\begin{document}

\title[Smooth Unbounded Worm]{Irregularity of the Bergman Projection
on Smooth Unbounded Worm Domains} 
\author[S. G. Krantz, A. Monguzzi, M. M. Peloso, C. Stoppato]{Steven G. Krantz, Alessandro Monguzzi, 
Marco M. Peloso, Caterina Stoppato}

\address{Campus Box 1146
Washington University in St. Louis
St. Louis, Missouri 63130}
\email{sk@math.wustl.edu}

 \address{Dipartimento di Ingegneria Gestionale, dell'Informazine e della Produzione
Università degli Studi di Bergamo
Viale G. Marconi 55
24044 Dalmine (BG) -- Italy}
\email{alessandro.monguzzi@unibg.it}

\address{Dipartimento di Matematica ``F. Enriques''\\
Universit\`a degli Studi di Milano\\
Via C. Saldini 50\\
I-20133 Milano}
\email{marco.peloso@unimi.it}

\address{Dipartimento di Matematica e Informatica ``U. Dini’’, Universit\`a degli Studi di Firenze, Viale Morgagni 67/A, I-50134 Firenze\bigskip}
\email{caterina.stoppato@unifi.it}


\thanks{{The second author is partially supported by the Hellenic Foundation for Research and Innovation (H.F.R.I.) under the “2nd Call for H.F.R.I. Research Projects to support Faculty Members \&
Researchers” (Project Number:  73342). He is also a member of  GNAMPA of the
  Istituto Nazionale di Alta Matematica {(INdAM)}.}}
\thanks{The  third author is partially
  supported by the GNAMPA Project CUP\_E55F2200027000.}
\thanks{{The fourth author is partly supported by GNSAGA of INdAM, by the INdAM project ``Hypercomplex function theory and applications'', by the PRIN project ``Real and Complex Manifolds'' and by the FOE project ``Splines for accUrate NumeRics: adaptIve models for Simulation Environments'' of the Italian Ministry of Education (MIUR)}}

 \begin{abstract}
In this work  we consider smooth unbounded worm domains $\cZ_\l$ in $\bbC^2$
and show that the Bergman projection, densely defined on the Sobolev
spaces $H^{s,p}(\cZ_\l)$,   $p\in(1,\infty)$,   $s\ge0$, does not extend to a
bounded operator $P_\l:H^{s,p}(\cZ_\l)\to H^{s,p}(\cZ_\l)$
when $s>0$ or $p\neq2$.  The same irregularity was known in the
case of the non-smooth unbounded worm. This improved result shows that
the irregularity of the projection is not a consequence of the
irregularity of the boundary but   instead of the infinite windings of the worm domain.
\end{abstract}
\maketitle

\section{Introduction}
Let $\phi$ be a non-negative smooth function on $\bbR$ such that
\begin{itemize}
   \itb $\phi$ is convex
   \itb $\phi^{-1}(0)=(-\infty, 0 ]$.
 %  ; \itb $\phi'<0$ on $(-\infty,-1)$.
\end{itemize}
Notice that $\phi'(t)>0$ for $t>0 $ and that
there exists $a> 0 $ such that $\phi(a)=1$.  For $\lambda>0$ we set
\begin{equation}\label{Zl-def}
\cZ_\l = \Big\{ (z_1,z_2):\,
\big|z_1-e^{i\log|z_2|^2}\big|^2<1-\phi(\log(\l|z_2|)^{-2}) \Big\}.
\end{equation}
Then, $\cZ_\l$ is smooth, unbounded  and pseudoconvex  (see Theorem 
\ref{main-thm1} below). Moreover, $\{\cZ_\l\}_{\l>0}$ is a nested
family of domains
whose union is the unbounded non-smooth worm
\begin{equation}\label{cW-def}
\cW = \Big\{ (z_1,z_2):\,
\big|z_1-e^{i\log|z_2|^2}\big|^2<1,\ z_2\neq0 \Big\}.
\end{equation}
The domain $\cW$ was studied in \cite{KPS1}, where three main facts were
proved   (see the enumerated list below). For $p\in[1,\infty]$ and $s\ge0$, given any domain $\Omega$,
denote by $H^{s,p}=H^{s,p}(\Omega)$ the
standard 
Sobolev space on $\Omega$. When $s=k$ is an integer, $H^{s,p}$ consists
  of functions with $k$-derivatives in $L^p(\Omega)$, and for
  non-integer $s$, $H^{s,p}$ can be defined by interpolation, see
  Section \ref{sec-2}.  For
    $p\in[1,\infty]$, 
let $A^p(\Omega):= L^p(\Omega)\cap\Hol(\Omega) $ denote the Bergman space.
    In
    \cite{KPS1} it was proved that:
    \begin{itemize}
    \item[(i)] the space $A^2(\cW)\neq\{0\}$, so that the Bergman
      projection $P: L^2(\cW)\to A^2(\cW)$ is a non-trivial orthogonal
      projector;
\item[(ii)] the operator $P$, initially defined on a dense subspace of $L^p(\cW)$, extends to a bounded operator $P:L^p(\cW)\to L^p(\cW)$ (if and) only if $p=2$; 
      \item[(iii)] the operator $P$, initially defined on a dense
subspace of $H^{s,2}(\cW)$, extends to a bounded operator
$P:H^{s,2}(\cW)\to H^{s,2}(\cW)$ (if and) only if $s=0$.
      \end{itemize}
The goal of this paper is to show that also in the case of the
unbounded smooth worms $\cZ_\l$, $\l>0$, the Bergman projection $P_\l$
on $\cZ_\l$ 
cannot be extended to a bounded operator $P_\l: H^{s,p}(\cZ_\l)\to
H^{s,p}(\cZ_\l)$ when $s>0$ or $p\neq2$. Observe that, since
$\cZ_\l \subseteq \cW$ for all $\l>0$, 
(i) above implies that $A^2(\cZ_\l)$ and hence $P_\l$ are
non-trivial.  We  now state our main results.

\begin{thm}\label{main-thm1}
  Let $\l>0$ and $\cZ_\l$ be defined as in \eqref{Zl-def}.  Then
$\cZ_\l$ is smooth, unbounded, and pseudoconvex and its boundary is
strongly pseudoconvex except at the points $\cA:= \{ (z_1,z_2): \,
z_1=0,\, |z_2|\ge 1/\l) \}$. Moreover, the
Bergman space $A^p(\cZ_\l)$ is infinite dimensional for all $p\in(0,\infty)$.
\end{thm}

\begin{thm}\label{main-thm2}
  Let $\l>0$, $\cZ_\l$ be defined as above, and let $P_\l$ denote the
  Bergman projection  on $\cZ_\l$.  If $P_\l$, initially defined on the
  dense subspace $(L^2\cap H^{s,p})(\cZ_\l)$, $p\in(1,\infty)$ and $s\geq0$,
extends to a bounded operator
  $$
  P_\l: H^{s,p}(\cZ_\l)\to
  H^{s,p}(\cZ_\l)
  $$
  then necessarily $s=0$ and $p=2$.   
\end{thm}

The problem of the regularity of the Bergman projection on worm
domains 
has been an object of active and intense research.  In the seminal paper
\cite{Barrett-Acta} D. Barrett considered the smoothly {\em bounded} worm domain
\begin{equation}\label{Wmu-def}
\cW_\mu=\left\{(z_1,z_2)\in\bbC^2: |z_1-e^{i\log
    |z_2|^2}|^2<1-\eta(\log|z_2|^2)\right\} ,
\end{equation}
 where
$\eta$ is smooth, non-negative, convex, 
$\eta^{-1}(0)=[-\mu, \mu]$,  and such that $\cW_\mu$ is smooth, bounded and pseudoconvex, see
e.g.~\cite[Proposition 2.1]{Ki}. Barrett showed that the
Bergman projection  on $\cW_\mu$ does not preserve the Sobolev space
$H^{s,2}(\cW_\mu)$ if $s\ge \pi/\mu$, whereas  in \cite{KP-Houston} it was
then shown that the
Bergman projection  on $\cW_\mu$ does not preserve $L^p$ if
$\big|\frac12-\frac1p\big|\ge\frac\pi\mu$. 
 We further mention in particular
 \cite{KP-Houston, MR2904008,  MR3424478, CS, KPS1, BDP}.
 We also refer the reader to \cite{MR2393268} for an 
 expository account of the subject, and to \cite{KPS2,PS1,PS2} for some
 interesting connections between Bergman spaces on worm domains and
 the {M}\"{u}ntz-{S}z\'{a}sz problem for the {B}ergman
              space in one complex dimension.

In the next section we prove Theorem \ref{main-thm1}, whereas in Section \ref{Sobolev:sec} we introduce the tools that we need to deal with Sobolev spaces on   smoothly bounded domains. In
Section~\ref{sec:3} we prove Theorem \ref{main-thm2} and in the final
Section~\ref{sec:last} we discuss some open problems and future work.
\ms

\section{The unbounded smooth worm}
\label{sec-2}

Consider the domains $\cZ_\l$.  
It is clear that they are unbounded, that
  $\cZ_\l\subseteq \cZ_{\l'}$ if $0<\l<\l'$ and that
  $\bigcup_{\l>0}\cZ_\l=\cW$. It is also immediate to see that
  $$
\cZ_\l \subseteq \{z_1:\, 0<|z_1|<2\} \times\{z_2:\, |z_2|>1/(e^{a/2}\l)\}. 
  $$
  Since $\cZ_\l\subseteq\cW$, where $\cW$
is as in \eqref{cW-def}, \cite[Proposition 2.3]{KPS1} gives that
$A^2(\cZ_\l)$ is infinite dimensional. Similar calculations also show
that also the spaces $A^p(\cZ_\l)$ are infinite dimensional,
$p\in(0,\infty]$.  Explicitly, for $\alpha\in\bbC$, for
$z=(z_z,z_2)\in\cW$, let
$$
L(z) = \log \big( z_1e^{-i\log |z_2|^2}\big) + i\log|z_2|^2,\qquad 
E_\alpha(z_1,z_2):= e^{\alpha L(z)},
$$
where $\log(z)$ denotes the principal branch of the logarithm on
$\bbC\setminus (-\infty,0]$. Then $L,E_\alpha\in\Hol(\cW)$ by
\cite[Lemma 2.2]{KPS1}.  Moreover,  for $j\in\bbZ$, $m\in\bbN$,
$c\in\bbR, c>\log2$, 
 $\alpha= \Re\alpha\CB +i(\frac{j}{2}+\frac1p)$, setting
$$
F_{\alpha,c,j,m}(z) := \frac{E_\alpha(z)z_2^j}{(L(z)-c)^m}
$$
and arguing as in  
\cite[Proposition 2.3]{KPS1}, it is simple to see that
$F_{\alpha,c,j,m}\in A^p(\cW)$ if $ \Re\alpha\CB >-2/p$ and $m>1/p$, where $p\in(0,\infty]$.  Hence,
$A^p(\cW)$ is infinite dimensional.

The argument to show that $\cZ_\l$ is smooth and pseudoconvex is
standard, but we repeat it for   the sake of completeness.  
Letting $\rho$
denote the defining function of $\cZ_\l$ we observe that
\begin{align*}
\rho(z_1,z_2) & = 
\big|z_1-e^{i\log|z_2|^2}\big|^2 -1+\phi(\log(\l|z_2|)^{-2})\\
& = |z_1|^2 -2\Re\big(z_1 e^{-i\log|z_2|^2}\big) +
\phi(\log(\l|z_2|)^{-2}),                                                             
\end{align*}
so that
$$
\p \rho(z_1,z_2)  =
\Big(\ov z_1 -e^{-i\log|z_2|^2}, -\frac{2}{z_2} \Im\big(z_1
e^{-i\log|z_2|^2}\big) -\frac{1}{z_2} \phi'(\log(\l|z_2|)^{-2}) \Big) .
$$
Let $(z_1,z_2)\in b\cZ_\l$ be such that $\p_{z_1} \rho(z_1,z_2)=0$.  Then $z_1
=e^{i\log|z_2|^2}$ so that $\phi(\log(\l|z_2|)^{-2})=1$.
The assumptions on $\phi$ imply that $\p_{z_2}\rho\neq0$ at such points.
Thus $\cZ_\l$ is smooth, and it is clearly unbounded since it
contains points $(z_1,z_2) $ with $|z_2|$ arbitrarily large.

  In order to show that $\cZ_\l$ is pseudoconvex, arguing as
  in~\cite{Ki}, we observe that  locally a branch of $\arg z_2$ is
  defined and that the
  local defining function $e^{\arg z_2^2}\rho$  \CB equals
$$
|z_1|^2 e^{\arg z_2^2} -2\Re \big(z_1 e^{-i\log z_2^2}\big)
+\phi(\log(\l |z_2|)^{-2}) e^{\arg z_2^2}.
$$
The first two terms are plurisubharmonic, while the third one
satisfies the differential inequality
$$
\Delta\big( \phi(\log(\l |z_2|)^{-2}) e^{\arg z_2^2} \big) =
\Delta\big( \phi(\log(\l |z_2|)^{-2})\big)  e^{\arg z_2^2} + \phi(\log(\l
|z_2|)^{-2}) \Delta\big( e^{\arg z_2^2} \big) \ge0,
$$
since $\phi$ is smooth and convex.  Hence, $\cZ_\l$ is pseudoconvex.
Moreover, the defining function is strictly plurisubharmonic at every
boundary point where $z_1\neq0$.

Next,  at $(z_1,z_2)\in b\cZ_\l$, the
complex tangent space  is spanned by the
vector
$$
 v= \begin{pmatrix} v_1 \\v_2 \end{pmatrix} := \begin{pmatrix}
 2\Im\big(z_1 e^{-i\log|z_2|^2}\big) +\phi'(\log(\l|z_2|)^{-2})
  \\
  z_2(\ov z_1 -e^{-i\log|z_2|^2}) \end{pmatrix}.
$$
Since
$$ \p^2_{z_1,\ov z_1} \rho =1,\
\p^2_{z_1,\ov z_2} \rho = \frac{i}{\ov z_2} e^{ - i\log|z_2|^2}, \ 
 \p^2_{z_2,\ov z_2} \rho = \frac{1}{|z_2|^2} \Big( 2\Re\big(z_1e^{-i\log|z_2|^2}\big) 
+ \phi''(\log(\l|z_2|)^{-2}) \Big), 
$$
and $2\Re(z_1 e^{-i\log|z_2|^2})= |z_1|^2 +
\phi(\log(\l|z_2|)^{-2})$ on the boundary,
the   
Levi form is given by
\begin{align*}
&   \cL_\rho(z; v)\\
  & =(v_1,v_2) \begin{pmatrix} 1 & \frac{i}{\ov z_2} e^{- i\log|z_2|^2} \\ -\frac{i}{z_2} e^{i\log|z_2|^2} & \frac{1}{|z_2|^2} \big( |z_1|^2 +
\phi(\log(\l|z_2|)^{-2})
+ \phi''(\log(\l|z_2|)^{-2}) \big) \end{pmatrix} \begin{pmatrix} \ov v_1\\
\ov v_2 \end{pmatrix} \\
& = v_1^2 + 2v_1 \Re\Big(i\frac{ e^{-i\log|z_2|^2}}{\ov z_2} \ov v_2 \Big)
 +|v_2|^2 \frac{1}{|z_2|^2} \Big( |z_1|^2 +
\phi(\log(\l|z_2|)^{-2})
+ \phi''(\log(\l|z_2|)^{-2}) \Big)\\
  & =  v_1^2  -2v_1 \Im\big(z_1e^{-i\log|z_2|^2} \big)
 +\big(1-\phi( \log(\l|z_2|)^{-2}) \big)\Big( |z_1|^2 +
\phi(\log(\l|z_2|)^{-2})
+ \phi''(\log(\l|z_2|)^{-2}) \Big)\\
& = 2\Im\big(z_1e^{-i\log|z_2|^2} \big) \phi'(\log(\l|z_2|)^{-2})
 + (\phi')^2 (\log(\l|z_2|)^{-2})\\
& \qquad\qquad +\big(1-\phi( \log(\l|z_2|)^{-2}) \big)\Big( |z_1|^2 +
\phi(\log(\l|z_2|)^{-2})
 + \phi''(\log(\l|z_2|)^{-2}) \Big).
\end{align*}
It follows that the boundary points $\{(0,z_2):
|z_2|\ge1/\l\}$ are of weak pseudoconvexity.
\begin{comment}
As in the case of the classical smooth worm, observe that
$(z_1,z_2)\in\cZ_\l$ if and only if
$|z_1e^{-i\log|z_2|^2}-1|^2<1-\phi( \log(\l|z_2|)^{-2})$, so that
$\Im(z_1e^{-i\log|z_2|^2})>0$.\footnote{\CR NO!} Thus,  the last term on the right in
the display above is non-negative, and therefore $\cZ_\l$
pseudoconvex.  Observe also that the boundary $b\cZ_\l$ is strongly
pseudoconvex except at that points $\{(0,z_2): |z_2|\ge 1/(\sqrt{e}\l)\}$. 
\end{comment}
This  proves Theorem \ref{main-thm1}.
\qed
\ms

\section{Sobolev spaces on   smoothly bounded domains and on $\cZ_\l$}\label{Sobolev:sec}
In this section we collect the results on Sobolev spaces on   smoothly
bounded domains and prove a few properties that we shall need later.
We begin by
recalling the definition and a few  standard
results from the theory  of function spaces on   smoothly {\em
  bounded} domains, see e.g.~\cite[Chapter
3]{Triebel} and~\cite{Lions-Magenes-artIII}.  In what follows the
space $H^{s,p}(\bbR^d)$
is defined by means of the
Fourier transform $\mathcal F$ on $\mathbb R^d$ and $\cD'(\Omega)$ is the dual of the space
$C^\infty_c(\Omega)$ of smooth functions with compact support in
$\Omega$. Namely,
\[
H^{s,p}(\bbR^d)=\Big\{f\in \cD'(\bbR^d): \mathcal F^{-1}\Big((1+|\xi|^2)^{\frac s2} \mathcal F f\Big)\in L^p(\bbR^d)\Big\}.
\]
\begin{defn}{\rm
Let $\Omega$ be a   smoothly bounded domain  in $\bbR^d$,  $s\ge0$ and
$p\in(1,\infty)$.  We define
$$
H^{s,p}(\Omega) = \Big\{ f\in \cD'(\Omega):\, \exists
F\in  H^{s,p}(\bbR^d)\,\big|\,  F_{|_\Omega}=f,\,
\|f\|_{H^{s,p}(\Omega)}:= \inf \{ \|F\|_{H^{s,p}(\bbR^d)}:\,
F_{|_\Omega}=f\} \Big\}.
$$
We also denote by $H^{s,p}_0(\Omega)$ the  closure of
$C^\infty_c(\Omega)$ in the $H^{s,p}(\Omega)$-norm.  Then, for $s<0$
and $p\in(1,\infty)$,
we define $H^{s,p}(\Omega)$ as the dual of $H^{-s,p'}_0(\Omega)$,
where $p'= p/(p - 1)$ is the
 exponent conjugate to $p$.\footnote{We remark that the definition of
   $H^{s,p}(\Omega)$ with $s<0$ is the same as in \cite{Lions-Magenes-artIII}
   but different from the one in \cite{Triebel}.}
}
\end{defn}

When $s=k$ is a non-negative integer the space $H^{k,p}(\Omega)$ has a
natural characterization. On the space $C^\infty(\ov\Omega)$ consider the norm 
$$
\| \psi\|_{W^{k,p}(\Omega)} := \sum_{|\alpha|\le k} \|\p^\alpha
\psi\|_{L^p(\Omega)} <\infty
$$
and define $W^{k,p}(\Omega)$ as the closure of 
$C^\infty(\ov\Omega)$ with respect to this norm.  Then
%% I eliminated the comma.
$W^{k,p}(\Omega)$ is isomorphic to $H^{k,p}(\Omega)$, with equivalence
of norms, see e.g.~\cite{Triebel}.

Using the complex interpolation method we have that, when $s>0$,
\begin{equation}\label{interpo:eq}
 H^{s,p}(\Omega) =[H^{k,p}(\Omega), H^{k+1,p}(\Omega)]_\theta 
\end{equation}
where $\theta\in(0,1)$ and $s=k+\theta$, cf. 
%% cfr. ---> cf.
\cite{Triebel} or~\cite{Lions-Magenes-artIII}, so that 
$[W^{k,p}(\Omega), W^{k+1,p}(\Omega)]_\theta$ is isomorphic as Banach
space to $H^{s,p}(\Omega)$, $s=k+\theta$.  For the complex
interpolation method we refer to \cite{Bergh-Lofstrom}.

\begin{comment}
Moreover, let $R$ denote the restriction operator $R: \cS'(\bbR^d)\ni
u\mapsto u_{|_\Omega}\in \cD'(\Omega)$, where $ \cS'(\bbR^d)$ denotes
the spaces of tempered distribution on $\bbR^d$. Then,
for every $N\in\bbN$, there exists a bounded linear operator 
$E$, called an {\em extension operator}, such that for all $0\ge s<N$ and $p\in(1,\infty)$ 
$$
E: H^{s,p}(\Omega)\to H^{s,p}(\bbR^d)
$$
such that
$RE=\operatorname{Id}$ on $H^{s,p}(\Omega)$.
\end{comment}

Since $\Omega$ is a bounded, smooth domain, the multiplier operator
$f\mapsto \chi_\Omega f$ is bounded on $H^{s,p}(\bbR^d)$ when
$0\le s<\frac1p$, $p\in(1,\infty)$. This fact in turn
%% turns ---> turn
implies the key property that $C^\infty_c(\Omega)$ is dense in
$H^{s,p}(\Omega)$ when $0\le s<\frac1p$  -- see~\cite[Theorem 3.4.3]{Triebel}.

We now prove a result that is probably well known, but for which we do
not know a precise reference.
\begin{lem}\label{mutual-duality:lem}
For $-1/p'< s<1/p$, the spaces $H^{s,p}(\Omega)$ and $H^{-s,p'}(\Omega)$
are mutually dual with respect to the $L^2(\Omega)$ pairing of duality.
\end{lem}  

\proof
Observe that, by duality, we may assume that $0\le s<1/p$. %%\footnote{Dettagli?}
Since $H^{s,p}(\Omega) =H^{s,p}_0(\Omega)$ in the given range,
$H^{-s,p'}(\Omega)=\big(H^{s,p}(\Omega))^*$
with the $L^2$-pairing of duality.

Conversely, let
$\ell\in\big(H^{-s,p'}(\Omega)\big)^* $.  Since the multiplication
$f\mapsto \chi_\Omega f$ is bounded on $H^{s,p}(\bbR^d)$,
$H^{s,p}(\Omega)$ can be identified with the subspace of
$H^{s,p}(\bbR^d)$ of functions vanishing on $\Omega^c$.
Therefore 
%% I eliminated the comma.
also $H^{-s,p'}(\Omega)$ can be identified with the
elements of $\big(H^{s,p}(\bbR^d)\big)^*=H^{-s,p'}(\bbR^d)$ that annihilate functions of 
$H^{s,p}(\bbR^d)$ vanishing on  $\Omega^c$.  Therefore, 
by the
Hahn--Banach theorem, there exists $L\in
\big(H^{-s,p'}(\bbR^d)\big)^*=H^{s,p}(\bbR^d)$
with the same norm, that agrees with $\ell$ on $H^{-s,p'}(\Omega)$.
Hence
%% I eliminated the comma.
 there exists $F\in H^{s,p}(\bbR^d)$ such that $\ell(u)=
\int_\Omega Fu = \int_\Omega (\chi_\Omega F)u$, where $\chi_\Omega
F\in H^{s,p}(\Omega)$; that is, $\big(H^{-s,p'}(\Omega)\big)^* =
H^{s,p}(\Omega)$.
\epf

Next we need an extension of a result by E. Ligocka, namely
\cite[Theorem 2]{Ligocka87}.
We denote by $H_{\har}^{s,p}(\Omega)$ the subspace of
$H^{s,p}(\Omega)$
consisting of harmonic functions.  Let $\varrho:\bbR^d\to\bbR$ be a smooth defining
function (see \cite{Kra1}) for $\Omega$ and let $L_{\har}^p(\Omega,|\varrho|^q)$ be the subspace
of $L^p(\Omega,|\varrho|^q \, dm)$ consisting of harmonic functions on
$\Omega$, $p\in(1,\infty)$.
In \cite[Theorem 2]{Ligocka87},
 Ligocka proved that,
   for $s\ge0$, $p\in(1,\infty)$,
\begin{itemize}
\item[(i)] $H_{\har}^{s,p}(\Omega)$ and
  $H_{\har}^{-s,p'}(\Omega)$ are mutually dual with respect to the
  $L^2(\Omega)$-inner product;
\item[(ii)] 
  $H_{\har}^{-s,p'}(\Omega)$  is isomorphically equivalent (as a Banach
  space) to $L_{\har}^{p'}(\Omega,|\varrho|^{sp'})$.
\end{itemize}

We shall need the following extension of (ii).

\begin{lem}\label{Ligocka:prop}
  Let $s\in\bbR$,
  $ s<1/p$, $p\in(1,\infty)$. Then,
  $H_{\har}^{s,p}(\Omega)$  is isomorphically equivalent (as a Banach
  space) to $L_{\har}^p(\Omega,|\varrho|^{-sp})$.
\end{lem}
\proof
As mentioned, the case $s\le0$ is proved in \cite[Theorem
2]{Ligocka87}. Next, let $0<s<1/p$.  
If $f\in L_{\har}^p(\Omega,|\varrho|^{-sp})$ and $g\in L_{\har}^{p'}(\Omega,|\varrho|^{sp'})$ we have that
\begin{align*}
\left |\int_\Omega fg\, dV \right |
  & = \left |\int_\Omega (|\varrho|^{-s} f)( |\varrho|^s g) \, dV\right |
    \le \|f\|_{L_{\har}^p(\Omega,|\varrho|^{-sp})}
    \|g\|_{L_{\har}^{p'}(\Omega,|\varrho|^{sp'})}\\
&    \le C \|f\|_{L_{\har}^p(\Omega,|\varrho|^{-sp})} \|g\|_{H_{\har}^{-s,p'}(\Omega)}
    ,
\end{align*}
so that, by (i) above,
\begin{align*}
 \| f\|_{H_{\har}^{s,p}(\Omega)} &\le C \| f\|_{(H_{\har}^{-s,p'}(\Omega))^*}
   = \sup\left \{ \left |\int_\Omega fg\, dV
\right |\,:\, \|g\|_{H_{\har}^{-s,p'}(\Omega)} \le 1\right \}\\
& \le C \|f\|_{L_{\har}^p(\Omega,|\varrho|^{-sp})} .
\end{align*}
Conversely, let $f\in H_{\har}^{s,p}(\Omega)$.  It is well known that the mapping
$H^{s,p}(\Omega)\ni f\mapsto |\varrho|^{-s}f \in L^p(\Omega)$ is
bounded when $0\le s<1/p$, see e.g.~\cite[Theorem 2, 1.3.1]{Mazya} or
\cite[p. 256]{Ligocka87}. Then we have,
\begin{align*}
\|f\|_{L_{\har}^p(\Omega,|\varrho|^{-sp})} & = \| |\varrho|^{-s}f \|_{ L^p(\Omega)}\\
& \le \|f\|_{H_{\har}^{s,p}(\Omega)}.
\end{align*} 
This proves the lemma.
\qed
\ms

We now define Sobolev spaces on the smooth unbounded domains $\cZ_\l$.
\begin{defn}{\rm
For  $k$ a non-negative integer and
$p\in(1,\infty)$, define the space (of test functions)
$$
\cT(\ov{\cZ_\l}) := \Big\{ \psi\in
C^\infty(\ov{\cZ_\l}):\, 
\| \psi\|_{H^{k,p}(\cZ_\l)} := \sum_{|\alpha|\le k} \|D_z^\alpha
\psi\|_{L^p(\cZ_\l)} <\infty \Big\},
$$
  where $D_z:=(\p_{z_1},\p_{\ov
  z_1};\p_{z_2}, \p_{\ov z_2})$. 
We define $H^{k,p}(\cZ_\l)$ as the closure of $\cT(\ov{\cZ_\l})$ with
respect to the norm $\| \cdot\|_{H^{k,p}(\cZ_\l)}$.  For $s=k+\theta$
with $0<\theta<1$, we define $H^{s,p}(\cZ_\l)$, $p\in(1,\infty)$, by
complex interpolation, as  
$$
H^{s,p}(\cZ_\l) := [H^{k,p}(\cZ_\l), H^{k+1,p}(\cZ_\l)]_\theta.
$$
See e.g.~\cite{Bergh-Lofstrom}. 
}
\end{defn}

Finally, we point out the following fact that we will need later.

\begin{remark}\label{ext-norm}{\rm
    Let $\mu(\l)=\log\l^2$, and consider the domain $\cW_{\mu(\l)}$ as defined in \eqref{Wmu-def},
where  $\eta$ is given by 
$$
  \eta(t)=\phi(t-\log\l^2)+\phi(-t-\log\l^2),
$$
so that  $\cW_{\mu(\l)}\subseteq\cZ_\l$. 
 Observe then that
the restriction operator $H^{s,p}(\cZ_\l)\ni f\mapsto
f_{|_{W_{\mu(\l)}}}\in H^{s,p}(\cW_{\mu(\l)})$ is well defined and norm
decreasing when $s=k$ is a non-negative integer and $p\in(1,\infty)$,
and then, by interpolation, also when $s\ge0$ and $p\in(1,\infty)$.
Analogously, for all $\l'>\l$,
$$
\| f\|_{H^{s,p}(\cW_{\mu(\l)})} \le \| f\|_{H^{s,p}(\cZ_{\l'})} .
$$
}
\end{remark}

\section{Irregularity of the Bergman projection}\label{sec:3}
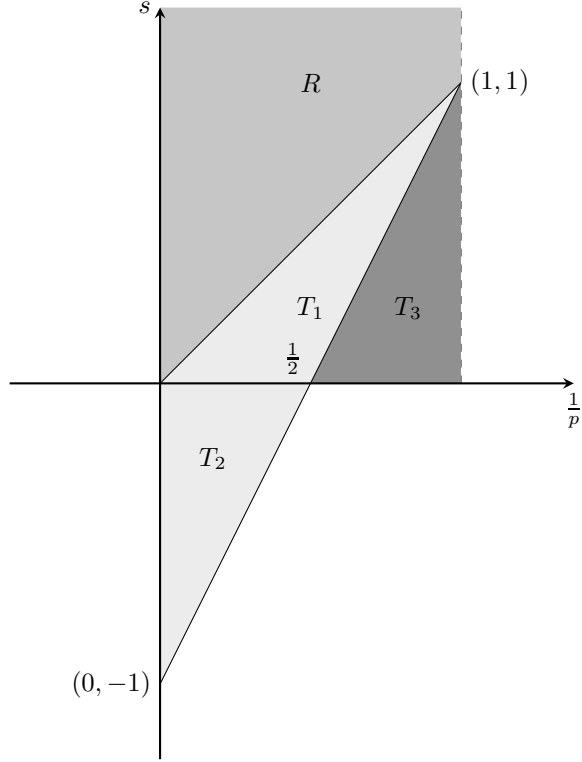
\begin{figure}
\begin{tikzpicture}[scale=1]
\draw[black, dashed] (4,0) -- (4,5);
 \fill[gray!45!white] (0,0) -- (0,5) -- (4,5) -- (4,4);
  \fill[gray!87!white] (2,0) -- (4,0) -- (4,4);
 \fill[gray!15!white] (0,-4) -- (0,0) -- (4,4) -- (2,0) -- cycle;
\draw[-stealth,black, thick] (-2,0) -- (5.5,0);
\draw[-stealth,black, thick] (0,-5) -- (0,5);
\draw[black] (0,0) -- (4,4);
\draw[black] (0,-4) -- (4,4);
\draw (5.5,0) node[below] {$\frac1p$};
\draw (0,5) node[left] {$s$};
\draw (4,4) node[right] {$(1,1)$};
\draw (2,0) node[above left] {$\frac12$};
\draw (0,-4) node[left] {$(0,-1)$};
\draw (2,4) node {$R$};
\draw (2,1) node {$T_1$};
\draw (0.7,-1) node{$T_2$};
\draw (3.3,1) node{$T_3$};
\end{tikzpicture}
\caption{Diagram for the proofs of Proposition \ref{extended-prop} and Theorem \ref{main-thm2}}
\label{fig1}
\end{figure}

 The proof of Theorem \ref{main-thm2} will combine some new ideas with
 Barrett's arguments \cite{Barrett-Acta} and results from
 \cite{KPS1}. We first extend \cite[Corollary 5.5]{KPS1} to the case
 of the Sobolev 
spaces $H^{s,p}(\cW_\mu)$.
\begin{prop}\label{extended-prop}
  Let $\cW$ be the unbounded non-smooth worm, $K_w$ its Bergman
  kernel at $w\in\cW$, and $\cW_\mu$ be the smoothly bounded worm as in \eqref{Wmu-def}.  Suppose
  $p\in(1,\infty)$.    Then the following properties hold:
  \begin{itemize}
  \item[(i)] if $s\in \big(
\frac2p -1,\infty\big)$ (the region $R\cup T_1\cup T_2$ 
  union the open
segments of end points $(0,0)$ and $(1,1)$ and $(0,0)$ and
$(\frac12,0)$, resp.,   in Figure
\ref{fig1}), then
 $K_w\not\in H^{s,p}(\cW_\mu)$;
    \item[(ii)] if $s=\frac2p -1$  and $p\in(1,2)$ (the open segment of end points
       $(\frac12,0)$, $(1,1)$ in Figure
\ref{fig1}),  then $\|K_w\|_{H^{s,p}(\cW_\mu)}\to
      \infty$ as $\mu\to\infty$. 
  \end{itemize}
\end{prop}

\proof
    We first observe that the cases 
$p=2$, 
$s>0$, and 
 $p>2$, $s=0$, that appear in (i),  are proved in \cite[Corollary
5.5]{KPS1}.

We now recall some notation from \cite[Corollary 5.5]{KPS1}.  We let
$S(e^{i\log|z_2|^2},\eps)$ denote the angular sector in the $z_1$-plane
$$
S(e^{i\log|z_2|^2},\eps) =\big\{ z_1=re^{i(t+\log|z_2|^2)}:\,
|t|<\delta,\, 0<r<\eps 
\big\} ,
$$
with $0<\delta< \pi/2$. 
For $\eps>0$ sufficiently small,
the set
\begin{equation}\label{gmu}
G_\mu =\big\{z= (z_1,z_2)\in\bbC^2: \,  |\log|z_2|^2| < \mu,  \, z_1\in  S(e^{i\log|z_2|^2},\eps)
\big\} 
\end{equation}
is contained in $\cW_\mu$. 
Then, from \cite[(5.8) and p. 1180]{KPS1},  for $w\in \cW$ and $z \in G_\mu$ we have the estimate
\begin{align}\label{p.1180}
  |K_w(z)|
  & \ge  \frac{C}{|z_1||z_2|} \CB \frac{1}{\big(\log(|z_1|/2)+\log(|w_1|/2)\big)^2 +  (\pi+2\mu)^2}
  \notag \\
  & \ge C_w\frac{1}{|z_1| |z_2|\CB ( \log^2 |z_1| +\mu^2)} ,
\end{align}
where $C_w$ does not depend on $\mu$.
Therefore, arguing as in \cite[Corollary 5.5]{KPS1}, for $s<1/p$ we have
\begin{align}
\| K_w \|^p_{L^p(\cW_\mu\CB,|\rho|^{-sp})}
  & \ge C_w \int_{|\log|z_2|^2|\le\mu}  \frac{1}{ |z_2|^p}\CB \int_{S(1,\eps)} 
\frac{1}{\big||\zeta|^2-2\Re\zeta\big|^{sp}  \big[ |\zeta|(\log^2|\zeta|
    +\mu^2)\big]^p} \,   dV(\zeta)\, dV(z_2) \notag \\
  & = C_w 2\pi   \frac{\sinh\big(|1-p/2| \mu\big)}{|1-p/2|} \CB \int_{|t|\le\delta}\int_0^\eps 
\frac{1}{|r-2\cos t|^{sp} r^{p(s+1)-1} (\log^2r+\mu^2)^p } \,
    drdt\notag \\
  & \ge C'_w    2\pi   \frac{\sinh\big(|1-p/2| \mu\big)}{|1-p/2|} \CB \int_0^\eps
    \frac{1}{r^{p(s+1)-1} (\log^2r+\mu^2)^p } \, dr. \label{est-Cor5.5}
\end{align}
\ms

 (i) Suppose then that 
$s\in\big(\frac2p -1,\frac1p\big)$. From Lemma
\ref{Ligocka:prop} $K_w\in H^{s,p}(\cW_\mu)$ if and only if $K_w\in L^p(|\varrho|^{-sp})$.
  From \eqref{est-Cor5.5} it then follows 
  that $K_w\not\in H^{s,p}(\cW_\mu)$  when $s\in(\frac2p -1,\frac 1p)$.  
We now use the natural embedding
$H^{s,p}(\cW_\mu)\subseteq H^{s',p}(\cW_\mu)$ when $0\le s' \le s$
(see~\cite[Theorem 3.3.1]{Triebel}). It follows that 
 $K_w\not\in H^{s,p}(\cW_\mu)$ for all
$p,s$ such that $p\in(1,\infty)$ and 
$s>\frac2p-1$. 
  This proves (i). 

 (ii) We look at the estimate
  in~\eqref{est-Cor5.5} when $s=\frac2p-1$ (notice that  $s<1/p$ in
this case) and observe that 
\begin{align}
\| K_w \|_{L^p(\cW_\mu,|\rho|^{-sp})}^p 
& \ge C'_w   2\pi   \frac{\sinh\big(|1-p/2| \mu\big)}{|1-p/2|}  \int_0^\eps
    \frac{1}{r (\log^2r+\mu^2)^p } \, dr \notag\\
& = C'_w   2\pi   \frac{\sinh\big(|1-p/2| \mu\big)}{|1-p/2|} 
\frac{1}{\mu^{2p-1}} \int_{\frac1\mu \log\frac1\eps}^\infty 
\frac{1}{(1+t^2)^p} \, dt \,. \label{est:Lp-weighted-bound}
\end{align} 
Clearly, if $p\neq2$, the right hand side above tends to $\infty$ if
$\mu\to\infty$. The rest of the proof will show that the same is true for $\|K_w\|_{H^{s,p}(\cW_\mu)}$. We observe in passing that, on the other hand,  if
$p=2$, the right hand side above
remains bounded (actually, it tends to $0$) when $\mu\to\infty$, in accordance to the
fact that $\|K_w\|_{L^2(\cW_\mu)}\le \|K_w\|_{L^2(\cW)}<\infty$. \ms

In order to conclude the proof of (ii), we will bound $\|K_w\|_{H^{s,p}(\cW_\mu)}$ from below. This will require several steps. We begin by setting, for $j=1,2$ and $\mu>2$,
$$
G_\mu^j =\big\{z= (z_1,z_2)\in\bbC^2: \, j< |z_2| < e^{\mu/2}/j,
 \, z_1\in  S(e^{i\log|z_2|^2},\eps/j)
\big\}\,.
$$
Keeping in mind~\eqref{gmu}, we see that $G_\mu^2\subseteq G_\mu^1\subseteq G_\mu$. 
We define a cut-off function $\psi:\bbC^2\to[0,+\infty)$ by setting $\psi(z_1,z_2) = \psi_1(z_1)\psi_2(z_2)$, where:
\begin{itemize}
\itb $\psi_1\in C^\infty_c(\bbC)$, $\psi_1(z_1) = 1$ for $|z_1|<
  \eps/2$ and $\psi_1 (z_1) = 0$ for $|z_1|\ge \eps$;
\itb  $\psi_2\in
C^\infty_c(\bbC)$ is identically equal to $1$ on
the annulus $\{z_2:\, 2<|z_2|<e^{\mu/2}/2\}$, is supported in a compact subset of
 $\{z_2:\, 1<|z_2|<e^{\mu/2}\}$ and has uniformly bounded derivatives.
\end{itemize}
%%Sobolev $p$-norms bounded independently of $\mu$.

Following the lines of the computations 
in~\eqref{est-Cor5.5} and~\eqref{est:Lp-weighted-bound}, we can find for each 
$p\in(1,2)$ and each $s=\frac2p-1>0$ a constant $C'>0$, independent of $\mu$, such that
$\| K_w\|_{L^p(G_\mu^2,|\rho|^{-ps})} \ge C' e^{\mu/p}$, whence
\begin{equation}\label{est:infty}
\| \psi K_w\|_{L^p(G_\mu^2,|\rho|^{-ps})} \ge C' e^{\mu/p}
\end{equation}
for $\mu>2$.  Now, let us consider the map $(z_1,z_2)\mapsto(z_1e^{-i\log|z_2|^2}, z_2)$. It is a $C^\infty$-diffeomorphism from (a neighborhood of) $G_\mu^1$ onto (a neighborhood of) its image
$$
\widetilde G_\mu^1 := \big\{\zeta= (\zeta_1,\zeta_2)\in\bbC^2: \, \zeta_1\in S(1,\eps),\, 1< |\zeta_2| <
e^{\mu/2} 
\big\}\,,
$$
included in the (Lipschitz) domain
$$
\Lambda =\big\{ w= (t+iu,w_2)\in\bbC^2:\, t>(\tan \delta)^{-1}|u|,\,
1< |w_2| < e^{\mu/2}  \big\}\,.
$$
We denote by $\Psi:\Lambda\to\bbC^2$ the inverse mapping $(w_1,w_2)\mapsto(w_1e^{i\log|w_2|^2}, w_2)$. We will later precompose $\Psi$ with the map $(w_1',w_2)\mapsto(\tau(w_1'),w_2)$ with $\tau(t'+iu)=t'+(\tan\delta)^{-1}|u|+iu$. The preimage of $\Lambda$ through this map (as well as $\Lambda$ itself) is contained in the half-space
$$
\cH:=\{ (\zeta,w_2)\in\bbC^2:\, \Re\zeta>0\}\,.
$$
We compute:
\begin{align*}
&  \| \psi K_w\|^p_{L^p(G_\mu^2,|\rho|^{-ps})}\le \| \psi K_w\|^p_{L^p(G_\mu^1,|\rho|^{-ps})}\\
&    =\int_{\widetilde G_\mu^1}  |(\psi K_w)\circ \Psi(w)|^p 
 |\rho\circ \Psi(w)|^{-ps}   |\det (J\Psi)(w)| \, dV(w)
    \\
& = \int_{\widetilde G_\mu^1}  |(\psi K_w)\circ \Psi(w) |^p
  ( 1-|w_1-1|^2)^{-ps}    |\det (J\Psi)(w)| \, dV(w)\\
 & \le \int_{\widetilde G_\mu^1}  |(\psi K_w)\circ \Psi(w) |^p
   ( 1-|w_1-1|^2)^{-ps}    \, dV(w)\\
& = \int_\Lambda | (\psi K_w)\circ \Psi(t+iu,w_2) |^p
 ( 2t-t^2-u^2)^{-ps}  \, dtdu \, dV(w_2)  \,,
\end{align*}
where we took into account the fact that $w\in\widetilde G_\mu^1$ implies $|w_2|>1$, whence $|\det (J\Psi)(w)|<1$.
On the support of $\psi\circ\Psi$, we have $t^2+u^2=|w_1|^2<\eps^2$. Up to shrinking $\eps$ to have $\eps<\big(1+(\tan\delta)^{-2}\big)^{-1}$, we get that $2t-t^2-u^2\ge t$ in the support of $\psi\circ\Psi$. We obtain
\begin{align*}
 \| \psi K_w\|^p_{L^p(G_\mu^2,|\rho|^{-ps})} 
  & \le \int_\Lambda | (\psi K_w)\circ \Psi(t+iu,w_2) |^p
 t^{-ps}  \, dtdu \, dV(w_2)\\
  & =\int_\cH | (\psi K_w)\circ \Psi(\tau(t'+iu),w_2) |^p
   ( t'+(\tan\delta)^{-1}|u|)^{-ps}  \, dt'du \, dV(w_2) \\
&   \le \int_\cH | (\psi K_w)\circ \Psi(\tau(t'+iu),w_2) |^p
   ( t')^{-ps}  \, dt'du \, dV(w_2) \\
&   = \| f\|^p_{L^p(\cH, (\Re \zeta)^{-ps})}\,,
\end{align*}
where $f(t'+iu,w_2) = (\psi K_w)\circ \Psi(\tau(t'+iu),w_2)$. For the half-space $\cH$, when $p\in(1,\infty)$, and $0<s<1/p$,
we have the well-known estimate
\begin{equation*}
\| h\|_{L^p(\cH, (\Re \zeta)^{-ps})} \le C'' \| h\|_{H^{s,p}(\cH)}
\end{equation*}
for all $h\in H^{s,p}(\cH)$, see e.g. \cite[Proposition 1 2.8.6,
Proposition 3.3.2]{Triebel}.  We conclude that 
\begin{align}\label{altra-stima}
 \| \psi K_w\|_{L^p(G_\mu^2,|\rho|^{-ps})}
\le C'' \| f\|_{H^{s,p}(\cH)} \,.
\end{align}
Our next aim is going back from $f(t'+iu,w_2) = (\psi K_w)\circ \Psi(\tau(t'+iu),w_2)$ to $(\psi K_w)\circ \Psi$. Taking into account that $\tau(t'+iu)=t'+(\tan\delta)^{-1}|u|+iu$, we compute
\begin{align*}
\p_u f(t'+iu,w_2) & = \p_u\Big( (\psi K_w)\circ
\Psi(\tau(t'+iu),w_2) \Big)\\
&  = (\tan\delta)^{-1}\sgn(u) \Big( \p_1 ((\psi K_w)\circ\Psi)\Big)(\tau(t'+iu),w_2)
% \\ & \qquad \qquad
                    +  \Big( \p_2 ((\psi K_w)\circ\Psi)\Big)(\tau(t'+iu),w_2)  \,,
\end{align*}
where we identified $\cH$ with $\bbR^4_+=\{(x_1,x_2,x_3,x_4):
\, x_1>0\}$ and we let $\p_1$ and $\p_2$ denote the partial derivatives w.r.t.~$x_1$ and $x_2$, resp. Hence,
\begin{align*}
& \big\| \p_u f(t'+iu,w_2) \big\|_{L^p(\cH)} \\
& \le C_\delta \Big( \Big\| \Big( \p_1 ((\psi K_w)\circ\Psi)\Big)(\tau(t'+iu),w_2) \Big\|_{L^p(\cH)} 
+ \Big\|  \Big( \p_2 ((\psi K_w)\circ\Psi)\Big)(\tau(t'+iu),w_2) \Big\|_{L^p(\cH)} \Big).
\end{align*}
Thus, there exists $C_\delta'>0$ such that, for all $g\in H^{1,p}(\cH)$ and for $\widetilde g(t'+iu,w_2)=g(\tau(t'+iu),w_2)$, the inequality
$$
\| \widetilde g\|_{H^{r,p}(\cH)} \le C_\delta' \|  g\|_{H^{r,p}(\cH)} \,.
$$
holds for $r=0,1$. By interpolation, the same inequality holds for all $r\in[0,1]$.
Using this bound in the estimate \eqref{altra-stima} and (later) the fact 
that all the derivates of the components of $\Psi$ are uniformly
bounded on the support of $\psi$, we obtain
\begin{equation*}
\| \psi K_w\|_{L^p(G_\mu^2,|\rho|^{-ps})} \le C \| (\psi K_w)\circ \Psi\|_{H^{s,p}(\cH)} \le C'\| \psi K_w\|_{H^{s,p}(\cW_\mu)}\,,
\end{equation*}
where all constants are independent of $\mu$. Now, the assumptions on $\psi_2$ guarantee that the function $\psi$ has $|D^\alpha_z\psi(z)|\le 1$ for all multiindices $\alpha$. Hence, multiplication by $\psi$ is a bounded operator, whose norm is independent of $\mu$, on $H^{k,p}(\cW_\mu)$ for all $k\in\bbN_0$ (whence on $H^{s,p}(\cW_\mu)$ for all
$s\ge0$). Thus, there exists a constant $C''$, independent of $\mu$, such that
\begin{equation*}
\| \psi K_w\|_{L^p(G_\mu^2,|\rho|^{-ps})} \le C'' \| K_w \|_{H^{s,p}(\cW_\mu)} \,.
\end{equation*}
This bound and~\eqref{est:infty} complete the proof of (ii).
\qed \ms

 In order to prove Theorem \ref{main-thm2} we need two preliminary lemmas. We denote by
 $\|T\|_{(X,X)}$ the operator norm of $T:X\to X$.

\begin{lem}\label{lem1}
For $\l, \l'>0$, the domain $\cZ_\l$ is biholomorphic to $\cZ_{\l'}$.  Moreover, the Bergman projection $P_\l$ 
induces a bounded
operator on $L^p(\cZ_\l)$ for some $\l>0$ if and only if $P_\l$ induces a bounded
operator on $L^p(\cZ_\l)$ for every $\l>0$ and in this case
$\|P_\l\|_{(L^p(\cZ_\l),L^p(\cZ_\l))}$ is independent of $\l$. 
\end{lem}
\proof
In order to show that the domains $\cZ_\l$ are all biholomorphic to each
other, it suffices to observe that for all $r,\l>0$, 
\begin{equation}\label{Phi-lambda-def}
\Phi_\l: \cZ_r\ni (w_1,w_2)\mapsto (w_1e^{-i\log\l^2}, w_2/\l) \in
\cZ_{\l r}
\end{equation}
is a biholomorphic map, since $\Phi_\l\in GL(2,\bbC)$ and
$\Phi_\l(\cZ_r)=\cZ_{\l r}$.  Moreover, $\det\Phi'_\l =
e^{-i\log\l^2}/\l$ and $T_{\l,p} f:=(\det\Phi'_\l )^{2/p} f\circ
\Phi_\l$ is an isometric isomorphism
$T_{\l,p}:  L^p(\cZ_{\l r})\to L^p(\cZ_r)$,
\begin{equation}\label{T-iso}
\| T_{\l,p} f\|_{L^p(\cZ_r)} = \|f\|_{L^p(\cZ_{\l r})}, 
\end{equation}
that also gives an   isometric isomorphism
$T_{\l,p}: A^p(\cZ_{\l r})\to A^p(\cZ_r)$ when restricted to
$A^p(\cZ_{\l r})$,
$p\in[1,\infty]$.  Recalling the transformation rule for the Bergman
projections
$$
P_r ( \det\Phi'_\l f\circ\Phi_\l) = \det \Phi'_\l (P_{\l r} f)\circ \Phi_\l
$$
for every $f\in L^2(\cZ_\l)$, 
since $\det \Phi'_\l$ is constant, it follows that $P_r (
f\circ\Phi_\l) =  (P_{\l r} f)\circ \Phi_\l$, for all $f\in L^2(\cZ_\l)$
and $\l>0$.  This implies that (also when $p\neq2$)
$$
P_r (T_{\l,p}f) =T_{\l,p}(P_{\l r} f)
$$
for all $f\in( L^2\cap L^p)(\cZ_{\l r})$.  Since $(L^2\cap
L^p)(\cZ_{\l r})$ is dense in $L^p(\cZ_{\l r})$ and $T_{\l,p}(L^2\cap
L^p)(\cZ_{\l r})$ is dense in $L^p(\cZ_r)$, for $f\in L^p(\cZ_{\l r})$ we have
$$
\| P_{\l r}f\|_{L^p(\cZ_{\l r})} = \| T_{\l,p}P_{\l r}f\|_{L^p(\cZ_r)} =  \|  P_r(T_{\l,p}f) \|_{L^p(\cZ_r)}. 
$$
Since $T_{\l,p}: L^p(\cZ_{\l r})\to L^p(\cZ_r)$ is an   isometric
isomorphism, the equality of the operator norms of $P_\l$ easily follows.
\qed\ms

\begin{lem}\label{lem2}
Let $s>0$, $p\in(1,\infty)$ and suppose $P_\l$ induces a bounded operator on
$H^{s,p}(\cZ_\l)$ for some $\l>0$.  Then, $P_{\l'}$ induces a bounded operator on
$H^{s,p}(\cZ_{\l'})$ for all $\l'>\l$ and, setting $N_{s,p}(\lambda)
=\|P_\lambda\|_{(H^{s,p}(\cZ_\lambda),H^{s,p}(\cZ_\lambda))}$, we have
\begin{equation}\label{ineq2}
N_{s,p}(\l')\le N_{s,p}(\l) .
\end{equation}
for all $\l'>\l$.  
\end{lem}

\proof
 For $r>0$, let
$T_{r,p}$ be as in the proof of Lemma \ref{lem1}. 
We argue as in~\cite{Barrett-Acta}. Recalling that $D_z=(\p_{z_1},\p_{\ov
  z_1};\p_{z_2}, \p_{\ov z_2})$,    if $\alpha=(a_1,b_1;a_2,b_2)$ is a
given multi-index, we have  that
$$
D_z^\alpha (f\circ \Phi_r) (z) = e^{i(b_1-a_1)\log r^2} r^{-(a_2+b_2)}
(D_z^\alpha f) (\Phi_r( z)) .
$$
Therefore, for $\l>0, r>1$,
 and $k$ a positive integer, using \eqref{T-iso}, we have
\begin{align}\label{ineq1}
  \| T_{r,p} f\|_{H^{k,p}(\cZ_\l)}
  & = \sum_{|\alpha|\le k} \| D_z^\alpha
T_{r,p} f\|_{L^p(\cZ_\l)} \le \sum_{|\alpha|\le k} \| 
T_{r,p} D_z^\alpha f\|_{L^p(\cZ_\l)}  =  \| 
f\|_{H^{k,p}(\cZ_{r\l})} .
\end{align}
Next observe that, using the transformation rule and a change of
variables, for $z\in\cZ_r$,
\begin{align*}
  (D_z^\alpha P_{r\l} f)(\Phi_r (z))
  & =  D_z^\alpha  \int_{\cZ_{r\l}} K_{r\l}(\Phi_r (z), w) f(w)\,   dV(w) \notag \\
  & = D_z^\alpha \int_{\cZ_{r\l}}
   |\det\Phi_r'|^{-2}   
    K_\l(z, \Phi_r^{-1} (w)) f(w)\,   dV(w) \notag\\
& =   D_z^\alpha \int_{\cZ_\l}
K_\l(z, w') f(\Phi_r(w'))\,   dV(w') \notag\\
  & = D_z^\alpha (P_\l(f\circ\Phi_r))(z) ,
\end{align*}
so that $T_{r,p} (D_z^\alpha P_{r\l} f)= D_z^\alpha(P_\l T_{r,p}f)$.

Therefore, assuming that $P_\l$ is bounded on $H^{s.p}(\cZ_\l)$, for
$r>1$, using both the fact that $T_{r,p}: L^p(\cZ_{r\l})\to L^p(\cZ_\l)$ is an isometry and \eqref{ineq1} we have
  \begin{align}  
    \| P_{r\l} f\|_{H^{k,p}(\cZ_{r\l)}}
    & = \sum_{|\alpha|\le k} \|D_z^\alpha P_{r\l} f\|_{L^p(\cZ_{r\l})} \notag\\
    & = \sum_{|\alpha|\le k} \|T_{r,p} D_z^\alpha P_{r\l} f\|_{L^p(\cZ_\l)}\notag\\
& =  \sum_{|\alpha|\le k} \| D_z^\alpha(P_\l T_{r,p}f)\|_{L^p(\cZ_\l)}  \notag\\
& = \|  (P_\l T_{r,p} f)\|_{H^{k,p}(\cZ_\l)} \notag\\
     & \le N_{k,p}(\l) \|   T_{r,p} f\|_{H^{k,p}(\cZ_r)} \notag\\
  & \le  N_{k,p}(\l) \|  f\|_{H^{k,p}(\cZ_{r\l})}          .            \label{id2} 
 \end{align}
Therefore $N_{k,p}(r\l)\le N_{k,p}(\l)$ for all integers $k$ and $r>1$. 
Thus, by interpolation, for all $s>0$ and $r>1$, \eqref{ineq2} follows.
\qed
\ms

\proof[Proof of Theorem \ref{main-thm2}]
{\em Step 1.} 
Suppose that $P_{\l_0}$ is bounded on $L^p(\cZ_{\l_0})$ for some
$\l_{0}>0$  and some $p\in(1,\infty)$. Hence $P_{\l}$ is bounded on
$L^p(\cZ_\l)$ for all $\l>\l_0$ 
by Lemma~\ref{lem1}.  Fix $f\in C^\infty_c(\cW)$ where $\cW$ is the
non-smooth unbounded worm and suppose that $\supp f\subseteq \cZ_\l$ for all
$\l\ge\l_0$.  For all such $\l$'s, denoting by $\chi_\l$ the
characteristic function of $\cZ_\l$,   
$$
\|\chi_\l P_\l f\|_{L^p(\cW)} = \| P_\l f\|_{L^p(\cZ_\l)} \le
C\|f\|_{L^p(\cZ_\l)} =C'
$$
for some constant $C'$ independent of $\l$. In the second-before-last
inequality, we have used Lemma \ref{lem1}. Then,  
 there exist a sequence $\{\l_n\}$, $\l_n\to\infty$ as
 $n\to\infty$, and $h\in (L^2\cap L^p)(\cW)$ such 
 that $\chi_{\l_n} P_{\l_n} f\to h$ in the weak-$*$ topology, as
 $n\to\infty$.  It is easy to see that $h\in
\Hol(\cW)$ arguing as follows. Let $\psi$ be  smooth, and compactly
supported in $\cW$. Then, denoting by $dV$ the Lebesgue volume
form,  for $j=1,2$ we have 
\begin{align*}
  \la (\ov\p_{z_j} h), \psi \ra
  & = 
\int_\cW (\ov\p_{z_j} h)  \ov\psi \, dV = -\int_\cW h
\ov{(\p_{z_j}\psi)}\, dV  =  
-\lim_{n\to\infty} \int_{\cZ_{\l_n}} P_{\l_n} f\ov{(\p_{z_j}\psi)}\, dV
  \\
  & =
\lim_{n\to\infty} \int_{\cZ_{\l_n}} (\ov\p_{z_j} P_{\l_n} f) \ov\psi \, dV= 
0.
\end{align*}
Hence, $\ov\p_{z_j} h=0$, $j=1,2$ and therefore $h$ is holomorphic.   We claim that
$h=Pf$, where $P$ denotes the Bergman projection on $\cW$.  It
suffices to show that $f-h\perp A^2(\cW)$. To this end, let $g\in
A^2(\cW)$. Then
$$
\int_\cW (f-h) \ov g \, dV  =  \lim_{n\to\infty} \int_\cW (f-\chi_{\l_n} P_{\l_n} f)
\ov g \, dV   = \lim_{n\to\infty} \int_{\cZ_{\l_n}} (f- P_{\l_n} f) \ov g \, dV    =0,
$$
since the restriction of $g$ to $\cZ_\l$ belongs to $A^2(\cZ_\l)$ for all $\l>0$ as well.

Seeking a contradiction,  we suppose $p\neq2$ and remark that
\begin{align*}
  \|P f\|_{L^p(\cW)}
  & = \sup \left \{ \left |\int_\cW Pf \ov\psi \, dV  \right | :\ \psi\in
    C_c^\infty(\cW),\|\psi\|_{L^{p'}(\cW)}\le1 \right \}\\
& =  \sup \left \{  \lim_{n\to\infty} \left |\int_{\cZ_{\l_n}} P_{\l_n} f \ov\psi \, dV  \right | :\ \psi\in
    C_c^\infty(\cW),\|\psi\|_{L^{p'}(\cW)}\le1 \right \}\\
  & \le  \lim_{n\to\infty}\|P_{\l_n} f\|_{L^p(\cZ_{\l_n})} \\
  & \le C \|f\|_{L^p(\cW)} .
\end{align*}
This implies that $P: L^p(\cW)\to L^p(\cW)$ is bounded, contradicting \cite[Theorem 1.1]{KPS1}.
Therefore, $P_{\l_0} $ cannot be bounded on $L^p(\cZ_{\l_0})$ for
$p\in(1,\infty)$ and $p\neq2$. We also observe that,
by interpolation with the
case $p=2$, $P_{\l_0}$ cannot be bounded on $L^1$ and $L^\infty$  either.
\ms

{\em Step 2.} In order to prove the irregularity of $P_\l$ in the Sobolev scale, 
we first show that $P_\l$ is densely
defined by showing that $(L^2\cap H^{s,p})(\cZ_\l)$ is dense in
$H^{s,p}(\cZ_\l)$. Let $\vp\in C^\infty_c(\bbC^2)$, $\vp=1$ on the
ball $B(0,1)$ and set $\vp^\eps ( \, \cdot \, ) =\vp(\eps\cdot)$. 
Given $f\in
H^{s,p}(\cZ_\l)$, let $f_{(\eps)}:=f\vp^\eps$.  It is easy to check
that $f_{(\eps)} \in (L^2\cap H^{s,p})(\cZ_\l)$ and that $f_{(\eps)}
\to f$ as $\eps \to 0^+$ in $H^{s,p}(\cZ_\l)$. \ms

{\em Step 3.} Let us show that it suffices to consider the case
$s\in(0,1/p)$ (the region $T_1\cup T_3$ in Figure \ref{fig1}). 
Suppose we have a bounded extension  $P_\l:H^{s,p}(\cZ_\l)\to H^{s,p}(\cZ_\l)$ for
some $s\ge 1/p$ and $p\in(1,\infty)$  (the region $R$ in Figure \ref{fig1}) . 
Interpolating with $L^2(\cZ_\l)$, we obtain a bounded extension $P_\l:H^{s_\theta,p_\theta}(\cZ_\l)\to H^{s_\theta,p_\theta}(\cZ_\l)$,
where $\theta\in(0,1)$, 
$s_\theta=\theta s$,
$\frac{1}{p_\theta}=\frac{\theta}{p}+\frac{1-\theta}{2}$. By taking
$\theta$ small enough we obtain that
$0<s_\theta<1/p_\theta$.
\ms

{\em Step 4.} We show that, 
if $p\in(1,\infty)$, $s\in(0,1/p)$, 
and $P_\l:H^{s,p}(\cZ_\l)\to H^{s,p}(\cZ_\l)$ is bounded, then $K_w\in
H^{s,p}(\cW_{\mu(\l)})$  and $K_w\in
H^{-s,p'}(\cW_{\mu(\l)})$. 

 Lemma \ref{lem2} gives bounded extensions
$P_{\l'}:H^{s,p}(\cZ_{\l'})\to H^{s,p}(\cZ_{\l'})$ for all $\l'>\l$
as well as 
$$
\|P_{\l'}\|_{(H^{s,p}(\cZ_{\l'}),H^{s,p}(\cZ_{\l'}))} \le N_{s,p}(\l)
$$
for all $\l'>\l$.

Fix $w\in\cW$ and let $K_w=K(\cdot,w)$ denote the Bergman kernel of $\cW$ at $w$. If we choose $\vp_w\in C^\infty_c$ supported in a ball
centered at $w$ within $\cW$, with radial symmetry and with $\int \vp_w =1$, then $P\vp_w = K_w$.

Then, for all $\l'>\l$ large enough for $\supp \vp_w\subseteq
\cZ_\l$, using Remark \ref{ext-norm} we have that
\begin{align*}
\| P_{\l'} \vp_w\|_{H^{s,p}(\cW_{\mu(\l)})}
  & \le \| P_{\l'} \vp_w\|_{H^{s,p}(\cZ_{\l'})} \le N_{s,p}(\l)
    \| \vp_w\|_{H^{s,p}(\cZ_{\l'})} =  N_{s,p}(\l)
    \| \vp_w\|_{H^{s,p}(\bbC^2)}.
\end{align*}  
Therefore $\{  P_{\l'} \vp_w\}_{\l'}$ is a family of functions contained in the
ball of radius $N_{s,p}(\l) \| \vp_w\|_{H^{s,p}(\bbC^2)}$ centered at the origin in
$H^{s,p}(\cW_{\mu(\l)})$. Since we are assuming $0<s<1/p$, 
using Lemma \ref{mutual-duality:lem} and the Hahn--Banach theorem we
have that $\{  P_{\l'} \vp_w\}_{\l'>\l}$ admits a subsequence
weak-$*$ converging to a function $h$ in $H^{s,p}(\cW_{\mu(\l)})$. Recalling
that $H^{s,p}(\cW_{\mu(\l)})$ is the dual of $H^{-s,p'}(\cW_{\mu(\l)})$ with respect to
 the
$L^2(\cW_{\mu(\l)})$ inner product, this implies that for all $g\in C^\infty_c(\cW_{\mu(\l)})$ we have
$$
\int_\cW (\chi_{\l'_n} P_{\l'_n}\vp_w)g \, dV  =
\int_{\cW_{\mu(\l)}}  (P_{\l'_n} \vp_w) g \to \int_{\cW_{\mu(\l)}}  h g \, dV  
$$
as $n\to\infty$.

Arguing as in {\em Step 1}, we have that (up to refinements) $\chi_{\l'_n} P_{\l'_n} \vp_w$ converges to $P \vp_w=K_w$ in the weak-$*$ 
topology of $L^2(\cW_{\mu((\l)})$. Thus
$$
\int_{\cW_{\mu(\l)}}  K_w g \, dV    = \int_{\cW_{\mu(\l)}}  h g \, dV  
$$
for all $g\in C^\infty_c(\cW_{\mu(\l)})$.  This implies that $h=K_w$ on
$\cW_{\mu(\l)}$, whence $K_w\in H^{s,p}(\cW_{\mu(\l)})$.

 In order to prove that $K_w\in
H^{-s,p'}(\cW_{\mu(\l)})$,
 we use Lemma \ref{mutual-duality:lem}.
For all $\l'>\l$ we have
\begin{align*}
\| P_{\l'} \vp_w\|_{H^{-s,p'}(\cW_{\mu(\l)})}
& = \sup \left \{ \left |\int_{\cW_{\mu(\l)}} P_{\l'} \vp_w\ov\psi \, dV  \right | :\ \psi\in
C_c^\infty(\cW_{\mu(\l)}), \, \|\psi\|_{H^{s,p}(\cW_{\mu(\l)})}\le1 \right \}\\
& =\sup \left \{ \left |\int_{\cZ_{\l'}} P_{\l'} \vp_w \ov\psi \, dV   \right | :\ \psi\in
    C_c^\infty(\cW_{\mu(\l)}),\, \|\psi\|_{H^{s,p}(\cW_{\mu(\l)})}\le1 \right \}  \\
  & = \sup \left \{ \left |\int_{\cZ_{\l'}} \vp_w \ov{P_{\l'} \psi} \, dV  \right | :\ \psi\in
    C_c^\infty(\cW_{\mu(\l)}),\, \|\psi\|_{H^{s,p}(\cW_{\mu(\l)})}\le1 \right \} \\
  & \le \| \vp_w\|_{H^{-s,p'}(\bbC^2)} \| P_{\l'} \psi\|_{H^{s,p}(\cZ_{\l'})} \\
  & \le N_{s,p}(\l) \| \vp_w\|_{H^{-s,p'}(\bbC^2)} .
\end{align*}
We now argue as before and conclude that $K_w\in
H^{-s,p'}(\cW_{\mu(\l)})$. 
\ms

We split the remaining part of the argument into 
three steps: one concerning the region $T_1$, one concerning the
region $T_3$ and one concerning the line segment separating them, see Figure \ref{fig1}.

{\em Step 5.} We assume $p\in(1,\infty)$ and $s>\max(\frac2p-1,0)$ and
$P_\l:H^{s,p}(\cZ_\l)\to H^{s,p}(\cZ_\l)$ is bounded. 
By {\em Step 4}, we have that $K_w\in H^{s,p}(\cW_{\mu(\l)})$.
Then, Proposition \ref{extended-prop} immediately gives a contradiction.
\ms

{\em Step 6.} We assume that $p\in(1,2)$,  
$0< s< \frac2p-1$ and $P_\l:H^{s,p}(\cZ_\l)\to
H^{s,p}(\cZ_\l)$ is bounded.
Notice that $\frac2p-1<\frac1p$, so that,
by {\em Step 4},
 we obtain that  $K_w\in H^{-s,p'}(\cW_{\mu(\l)})$, where
$-s > 1-\frac2p = \frac{2}{p'}-1 $.
But, again, this is false by Proposition \ref{extended-prop} 
(i)   and we
have reached a contradiction. Hence,
 the projector $P_\l$ does not
extend to  a bounded operator $P_\l :H^{s,p}(\cZ_\l)\to H^{s,p}(\cZ_\l)$.

  {\em Step 7.} Finally, let $p\in(1,2)$ and $s=\frac2p-1$ and suppose
that $P_\l:H^{s,p}(\cZ_\l)\to
H^{s,p}(\cZ_\l)$ is bounded.  Again, Lemma \ref{lem2} gives that
$P_{\l'}:H^{s,p}(\cZ_{\l'})\to H^{s,p}(\cZ_{\l'})$ is bounded 
and 
$$
\|P_{\l'}\|_{(H^{s,p}(\cZ_{\l'}),H^{s,p}(\cZ_{\l'}))} \le N_{s,p}(\l)
$$
for all $\l'>\l$.  Take any $\mu>\pi$ and let $\l'$ sufficiently large so
that $\cW_\mu\subseteq \cZ_{\l'}$. 
Let $\vp_w\in C^\infty_c(\cZ_{\l'})$ be as in {\em Step 4}. We have shown that there
exists a sequence $\{ P_{\l_n}\vp_w\}$ such that $P_{\l_n}\vp_w\to
K_w$ weak-$*$ in $H^{s,p}(\cZ_{\l'})$ as $\l_n\to\infty$, so that
$$
\| K_w\|_{H^{s,p}(\cW_\mu)} \le 
\| K_w\|_{H^{s,p}(\cZ_{\l'})} \le \lim_{n\to\infty}
\|P_{\l_n}\|_{(H^{s,p}(\cZ_{\l_n}),H^{s,p}(\cZ_{\l_n}))}
\|\vp_w\|_{H^{s,p}(\bbC^2)} \le C  N_{s,p}(\l)
$$
independent of  $\mu$.  This contradicts Proposition \ref{extended-prop} (ii) and the proof is complete.
\qed
\ms

\section{Final remarks and open questions}
\label{sec:last}  

We wish to conclude by indicating a number of open problems.  First of
all, we recall that the exact range of regularity on the Lebesgue
Sobolev spaces $H^{s,p}$ of the Bergman
projection on the   smoothly bounded domain $\cW_\mu$ is not known.
Clearly, in order to prove a positive result, one needs to have
precise information on the Bergman kernel itself.  In fact, also the precise
behaviour of the kernel near the critical annulus $\cA=\{ (0,z_2):\,
e^{-\mu/2}<|z_2|<e^{\mu/2}\}$ on the boundary of $\cW_\mu$ remains to be
understood.

The equivalence of the regularity of the Bergman projections on
$(0,q)$-forms and the Neumann operator $\cN$, proved in
\cite{BS-equivalence}, was later exploited by M. Christ \cite{Christ}
to show that $P_\mu$ does not preseve $C^\infty(\ov{\cW_\mu})$.  These
results heavily relied on the boundedness of the domain $\cW_\mu$. We
believe that the Neumann operator $\cN$ on $\cZ_\l$ is as irregular as
the Bergman projection $P_\l$, but this problem has not been addressed
and (to the best of our knowledge) is open.

Finally, we mention the boundary analogue of this problem, namely the
study of the behaviour of the Szeg\H{o} projection on $\cZ_\l$.
Given a smooth domain $\Omega=\{ z:\rho(z)<0\}
\subseteq\bbC^n$, the Hardy space $H^2(\Omega,d\sigma)$ is 
defined as 
$$
H^2(\Omega,d\sigma) =\big\{ f\in\Hol(\Omega): \sup_{\eps>0}
\int_{\p\Omega_\eps} |f|^2 d\sigma_\eps <\infty
\big\} \,,
$$
where $\Omega_\eps=\{ z:\rho(z)<-\eps\}$ and $d\sigma_\eps$ is the
induced surface measure on $\p\Omega_\eps$. 
Then $H^2(\Omega,d\sigma)$ can be 
identified with a closed subspace of $L^2(\p\Omega,d\sigma)$, that we
denote by $H^2(\p\Omega,d\sigma)$, where $\sigma$ is the induced
surface measure on $\p\Omega$.  The
Szeg\H o projection is the orthogonal projection 
$$
S_\Omega : L^2(\p\Omega,d\sigma) \to H^2(\p\Omega,d\sigma) \,;
$$
see \cite{St2} for the case of {\em bounded} domains.  The regularity
of $S_\Omega$ when $\Omega$ is a (model) worm domain was studied in a series
of papers \cite{MR3506304,MR3491881,M,MP,MP2,LS-Cortona}. In
particular, in  \cite{LS-Cortona} it was announced that $S_{\cW_\mu}$
does not preserve $L^p(\p\cW_\mu)$ when
$\big|\frac12-\frac1p\big|\ge\frac\pi\mu$, in analogy to the case of the
Bergman projection.  L. Lanzani and E. Stein also studied the
$L^p$-regularity of the Szeg\H{o} and other projections on the
boundary on bounded domains under minimal smoothness conditions
\cite{MR3250300,2015arXiv150603965L}, whereas a definition of Hardy
spaces and associated Szeg\H o projection for singular domains was
studied, for instance, in \cite{monguzzi-Hartogs, GGLV}. It is
certaintly of interest to consider the case of the Szeg\H o 
projection also in the case of the domains $\cZ_\l$. 
\ms

\subsection*{Declarations}
Data sharing not applicable to this article as no datasets were generated or analysed during the current study.
The authors have no relevant financial or non-financial interests to disclose.

\bibliography{worm-mathscinet}
\bibliographystyle{amsalpha}

\end{document}